\documentclass[twoside,11pt]{article}
\usepackage{graphicx}
\usepackage{amssymb}
\usepackage{amsthm}
\newfont{\bbb}{msbm10 scaled\magstep 1}
\newcommand{\bd}{{\rm bd}}

\newcommand{\conv}{{\rm conv}}
\newcommand{\width}{{\rm width}}
\newtheorem{thm}{Theorem}

\font\bigbold=cmbx10 at 14 pt
 
\date{}

\begin{document}

\title{\bf Banach-Mazur distance from the parallelogram to the affine-regular hexagon and other affine-regular even-gons}
 
\baselineskip 17.5 pt 

\maketitle

\vskip -1.3cm 
\centerline
{\bigbold Marek Lassak}

\vskip0.25cm
\centerline {\it University of Science and Technology}

\centerline {\it Kaliskiego 7, 85-789 Bydgoszcz, Poland}

\centerline {\it e-mail: lassak@utp.edu.pl}

\vskip0.65cm

\pagestyle{myheadings} \markboth{ Marek Lassak}{Banach-Mazur distance}

\vskip 0.3cm 

\begin{abstract}
We show that the Banach-Mazur distance between the parallelogram and the affine-regular hexagon is $\frac{3}{2}$ and we 
conclude that the diameter of the family of centrally-symmetric planar convex bodies is just $\frac{3}{2}$.
A proof of this fact does not seem to be published earlier. 
Asplund announced this without a proof in his paper proving that the Banach-Mazur distance of any planar
centrally-symmetric bodies is at most $\frac{3}{2}$. 
Analogously, we deal with the Banach-Mazur distances between the parallelogram and the remaining affine-regular even-gons.
\end{abstract}

\section {Introduction}
\label{Intro}

Denote by $\mathcal{M}^d$ the family of all centrally symmetric convex bodies of the Euclidean space $E^d$ centered at the origin $o$ of $E^d$.
Below always when we say on a homothety, we mean that its center is $o$.
For any $K \in \mathcal{M}^d$ and any positive $\lambda$, by $\lambda K$ we denote the homothety image of $K$ with ratio $\lambda$.

The Banach-Mazur distance (or shortly, the BM-distance) of $C, D \in \mathcal{M}^d$ is defined as 

$$\delta_{BM} (C,D) = \inf_{a, \lambda} \{ \lambda ; \, a(C) \subset D \subset \lambda a(C) \}, $$

\noindent
where $a$ stands for an affine transformation. 
This definition is presented by Banach \cite{B} in behalf of him and Mazur. 
From the role of the affine transformation in this definition we see that a different formulation of this definition is formally more precise, namely when we consider the Banach-Mazur distance of the equivalence classes of centrally symmetric bodies with respect to affine transformations.
But in this note it is more convenient to deal with the BM-distance of convex bodies from $\mathcal M^d$.
It is well known that the Banach-Mazur distance is a multiplicative metric, i.e., that ${\rm log} \, \delta_{BM}$ is a metric.
In particular, we have $\delta_{BM} (C,D) = \delta_{BM} (D,C)$ and the multiplicative triangle inequality $\delta_{BM} (C,D) \hskip0.02cm  \cdot \hskip0.02cm  \delta_{BM} (D,E) \geq \delta_{BM} (C,E)$ for every $C, D, E \in \mathcal{M}^d$.  
Later the notion of the BM-distance has been generalized to pairs of arbitrary convex bodies.
For a survey of results on the BM-distance we propose
the books by Tomczak-Jaegerman \cite{TJ}, Toth \cite{T} (Part 3.2), and Aubrun and Szarek \cite{ASz} (Chapter 4).

By $P_n$ we  denote the regular $n$-gon with vertices $(\cos \frac{2j}{n}\pi,  \sin \frac{2j}{n}\pi)$ for $j= 0, \dots , n-1$.

Our main aim is to show that $\delta_{BM} (P_4, P_6) = \frac{3}{2}$ and to establish all the optimal positions of the parallelograms $a(P_4)$ with respect to the hexagon $P_6$. 
Recall that the value $\frac{3}{2}$ is stated in Theorem 1 of \cite{A}.
It says that $\delta_{BM} (P_4, C) \leq \frac{3}{2}$ for any $C \in \mathcal M^2$ with equality only for $C= P_6$.
But the proof of Lemma 3 of \cite{A} by Asplund that $\delta_{BM} (P_4, P_6) = \frac{3}{2}$ is omitted.
This statement is quoted on p. 206 by Stromquist \cite{S} and illustrated in Fig. 1 there.
Since the hexagon from this Fig. 1 is not inscribed in the larger parallelogram, then looking to the second thesis of our Proposition we get further questions. 
The basic task of \cite{S} is to construct a center $R$ of $\mathcal M^2$ (see p. 207 and compare Fig. 2). 
With the help of Fig. 3 it is shown that both $P_4$ and $P_6$ are in the distance $\sqrt {3/2}$ from $R$. 
This implies $\delta_{BM} (P_4, P_6) \leq \frac{3}{2}$, but does not imply the equality.
These doubts mobilized the writer to present a detailed proof of the equality $\delta_{BM} (P_4, P_6) = \frac{3}{2}$ in the present note in order to be sure that the diameter of $\mathcal M^2$ is~$\frac{3}{2}$.

An additional aim is explained here.
The idea of the proof of our Theorem \ref{hexagon} that $\delta_{BM} (P_4, P_6) = \frac{3}{2}$ encouraged the author to consider the analogous task for all the regular even-gons with more vertices in place of $P_6$.
Theorem 2 establishes the Banach-Mazur distances from $P_4$ to all $P_{8j}$ and $P_{8j+4}$. 
This task for the remaining odd-gons, so to $P_{8j+2}$ and $P_{8j+6}$ appears to be more complicated.
We estimate and conjecture the values of these distances.

\section{Positions of parallelograms which realize the Banach-Mazur distance to a given centrally symmetric convex body}

By an inscribed parallelogram in a convex body $C$ we mean a parallelogram with all vertices in the boundary of $C$ and by a circumscribed parallelogram about $C$ we mean a parallelogram containing $C$ whose all sides have non-empty intersections with $C$.

By the width $\width(S)$ of a strip $S$ between two parallel hyperplanes of $E^d$ we mean the distance of these hyperplanes. 
We omit an easy proof of the following lemma, whose two-dimensional version is applied in the proofs of Theorems \ref{hexagon} and \ref{I-IV}.

\vskip0.25cm
\noindent
{\bf Lemma.}
{\it Let $H^1_+, H^1_-$ and $H^2_+, H^2_-$ be two pairs of hyperplanes of $E^d$ symmetric with respect to $o$ such that all of them are parallel and let $L$ be a straight line through $o$ intersecting them. 
Assume that $H^2_+, H^1_+,  H^1_-, H^2_-$ intersect $L$ in this order.
Denote by $(a^i_1, \dots , a^i_n)$ the points of intersection of $L$ with $H^i_+$, where $i\in \{1, 2\}$.
Consider the strips $S^i = \conv (H^i_+ \cup H^i_-)$, where $i\in \{1, 2\}$.
We claim that $\width(S^2)/\width(S^1) = a^2_j/a^1_j$ for any $j \in \{1, \dots , n\}$.} 

\medskip
The following proposition is applied later for evaluating particular BM-distances.

\medskip
\noindent
{\bf Proposition.}
{\it Let $C \in \mathcal{M}^2$.
Assume that $P \subset C \subset \mu P$ for a parallelogram $P \in \mathcal{M}^2$ and a positive $\mu$. 
Then there exists a parallelogram $P'$ inscribed in $C$ such that $\mu'P'$ is circumscribed about $C$ for a
$\mu' \leq \mu$.}

\begin{proof} 
Since the thesis is obvious for $\mu =1$, below we assume that $\mu > 1$.

If $P$ in the part of $P'$ fulfills the thesis, there is nothing to prove.
In the opposite case we intend to construct the required $P'$.

Of course, there is a homothetic image $P_\alpha \subset P$ of $P$ 
such that for least one pair of opposite sides $Z^+, Z^-$ of $P_\alpha$
the sides $\mu Z^+, \mu Z^-$ of $\mu P_\alpha$ touch $C$ from outside.
Denote the other pair of sides of $P_\alpha$
by $W^+, W^-$ and assume that their order is $W^+, Z^+, W^-, Z^-$when we move counterclockwise.

If $\mu W^-, \mu W^+$ do not touch $C$, then we lessen $P_\alpha$ up to $P_\beta$ by moving $W^+, W^-$
such that their vertices remain in $Z^+, Z^-$ 
symmetrically closer to $o$ up to the position $F^+ , F^-$ when both $\mu F^+$ and $\mu F^-$ touch $C$.
The other two sides of the parallelogram $P_\beta$ are denoted by $G^+, G^-$.
Clearly, $G^+ \subset Z^+$ and $G^- \subset Z^-$.
The parallelogram $\mu P_\beta$ is circumscribed about $C$.
Of course, when we go counterclockwise, then the order of the sides of $P_\beta$ is $F^+, G^+, F^-, G^-$.

Clearly $P_\beta$ is a subset of $C$, but it may be not inscribed in $C$.

Take the homothetic copy $P_\gamma \subset C$ of $P_\beta$ such that at least one pair of opposite vertices $t, v$ of $P_\gamma$ is in $\bd (C)$.
Denote the other two vertices of $P_\gamma$ by $u, w$ such that $t, u, v, w$ be in this order on $\bd (C)$.  

Let 
$a, b \in \bd (C)$ be the points such that $w \in ta$ and $w \in vb$

If $u, w$ are not in $\bd (C)$, let us cleverly enlarge $P_\gamma$ up to a parallelogram $P_\delta$ (not obligatory homothetic) inscribed in $C$ such that a homothetic copy of $P_\delta$ with a ratio at most $\mu$ is circumscribed about $C$.

Here we explain how to provide this task.
For every boundary point $c$ of $\bd (C)$ on the arc $\stackrel{\frown}{ab}$ provide the straight lines 
$K_1^+(c), K_2^+(c)$ containing $cv$ and $ct$, respectively.
Let $K_1^-(c), K_2^-(c)$ be the line symmetric to $K_1^+(c), K_2^+(c)$, respectively.
Moreover, for $i=1, 2$ provide the two pairs of the straight lines $L_i^+(c)$ and $L_i^-(c)$ supporting $C$ and parallel to the pairs $K_i^+(c)$ and $K_i^-(c)$ such that the order of them is $L_i^+(c)$, $K_i^+(c)$, $K_i^-(c)$, $L_i^-(c)$.

Let $L_i(c)$ be the strip between $L_i^+(c)$ and $L_i^-(c)$ for $i= 1, 2$
and $K_i(c)$ be the strip between $K_i^+(c)$ and $K_i^-(c)$ for $i= 1, 2$.

Provide the straight line through $o$ and points $d^+, d^-$ of support of $C$ by $L_1^+(c)$ and $L_1^-(c)$, respectively.
Denote by $f^+, f^-$ its intersections with $F^+, F^-$, respectively.
Denote by $k^+, k^-$ its intersections with $K_i^+(c), K_i^-(c)$, respectively.
Denote by $h^+, h^-$ its intersections with $\mu F^+, \mu F^-$.
 Since $h^+$, $d^+$, $k^+$, $f+$, $o$, $f^-$, $k^-$, $d^-$, $h^-$ are in this order on this straight line, we obtain $|h^+h^-|/|f^+f^-| \leq |d^+d^-|/|k^+k^-| $.
Hence by Lemma we obtain that 

$$\width (L_i(c))/\width (K_i(c)) \leq \mu$$ 

\noindent
for $i=1$.
Analogously, we conclude this for $i=2$.

If $c$ is sufficiently close to $a$, then 

$$\width (L_1(c))/\width (K_1(c)) > \width (L_2(c))/\width (K_2(c)).$$

If $c$ is sufficiently close to $b$, then we have the opposite inequality.
Moreover, observe that when $c$ moves on $\stackrel{\frown}{ab}$ from $a$ to $b$, then the strips $L_i(c)$ and $K_i(c)$, where $i = 1, 2$, change continuously. 
Consequently, there is at least one position $c_0$ of $c$ for which 

$$\width (L_1(c_0))/\width (K_1(c_0)) = \width (L_2(c_0))/\width (K_2(c_0)).$$

Therefore the thesis of our proposition is true for the parallelograms $P' = K_1(c_0) \cap K_2(c_0)$. 
Still the homothetic one $\mu'P = L_1(c_0) \cap L_2(c_0)$, where $\mu' \leq \mu$, is circumscribed about $C$.
\end{proof}

\medskip
\noindent
{\bf Corollary.}
{\it Denote the Banach-Mazur distance between $C \in \mathcal{M}^2$ and $P$ by $\mu$. 
Assume that $P \subset C \subset \mu P$ for a particular affine image $P \in \mathcal{M}^2$ of $P_4$. 
Then the parallelogram $P$ is inscribed in $C$ and $\mu P$ is circumscribed about $C$.} 

\medskip
The author does not know if Proposition holds true in higher dimensions for the parallelotope or the cross-polytope in place of the parallelogram.

\section{Banach-Mazur distance between the parallelogram and 
the affine-regular hexagon}

\begin{thm}\label{hexagon}
We have $\delta_{BM} (P_4, P_6) = \frac{3}{2}$. 
\end{thm}

\begin{proof}
The parallelogram with the four vertices at $(\pm 1, 0)$ and $(0, \pm \frac{1}{2}\sqrt 3)$ is contained in $P_6$ and its homothetic image with ratio $\frac{3}{2}$ contains $P_6$.
Consequently, $\delta_{BM} (P_4, P_6) \leq  \frac{3}{2}$.

Having in mind Proposition, in order to show that $\delta_{BM} (P_4, P_6) \geq  \frac{3}{2}$, it is sufficient to consider only any parallelogram $P = a(P_4)$ inscribed in $P_6$ such that a positive homothetic copy $\lambda P$ of $P$ is circumscribed about $P_6$, and to show that this homothety ratio $\lambda$ is at least $\frac{3}{2}$.  

Denote by $\mathcal{P}$ the class of all such parallelograms $P$. 

Consider a parallelogram $P = pqrs$ from $\mathcal{P}$. 
Some two consecutive vertices of $P$ must be in two consecutive sides of $P_6$.
The reason is that in the opposite case no positive homothetic image of our $P$ is circumscribed about $P_6$ in contradiction to $P \in \mathcal{P}$.
In order to fix attention, thanks to the symmetries of $P_6$, we do not make our considerations narrower
assuming that $p \in v_0v_1$ and $q \in v_1v_2$.
For the same reason, we may additionally assume that $p \in v_0m$, where $m$ denotes  the middle of the side $v_0v_1$.

{\ }

\begin{center}

\includegraphics[width=3.8in]{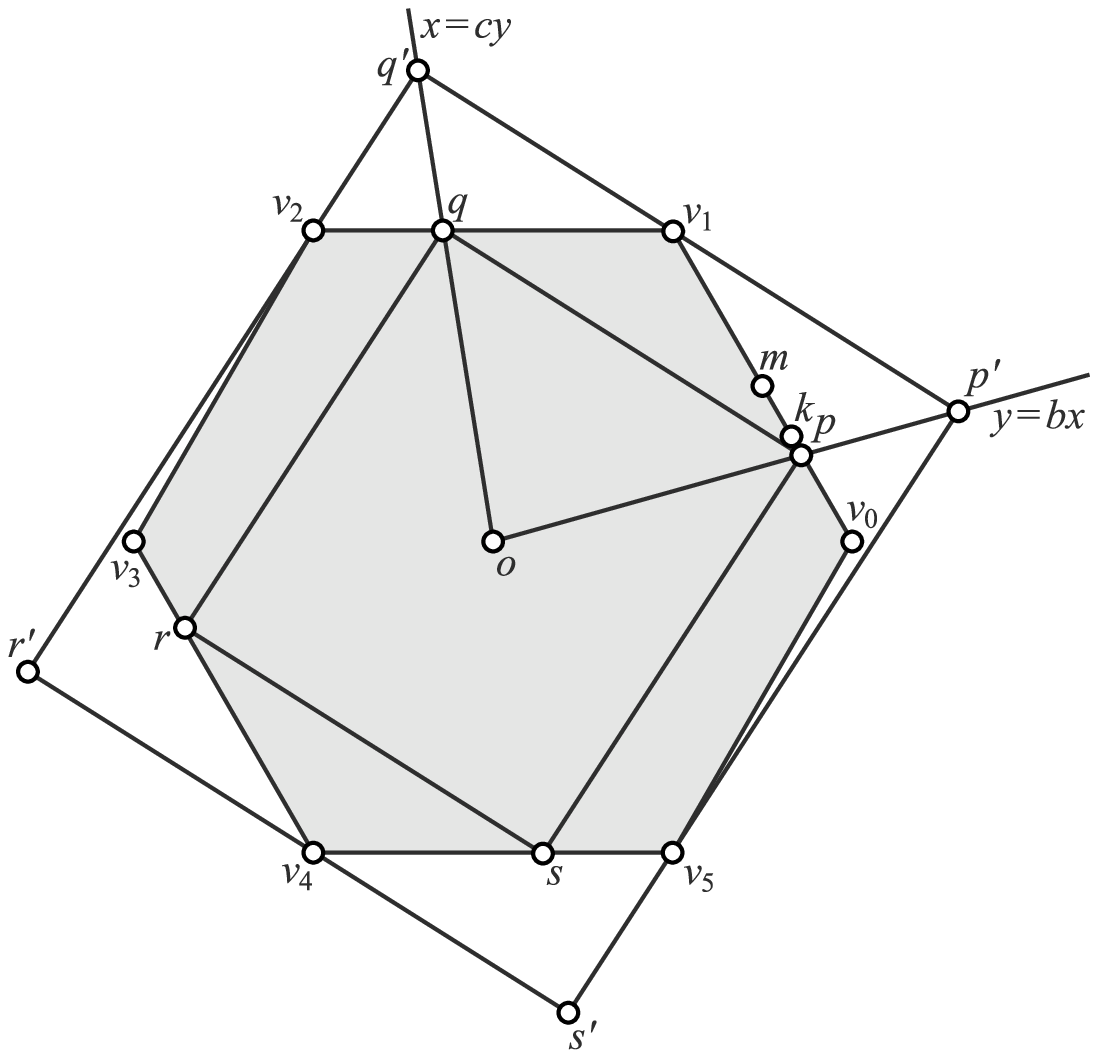} \\ 

\vskip0.2cm
{Fig 1. The positions of $P \in \mathcal{P}$ and $h(a)P$ with respect to the hexagon} 

\end{center}

Take the line $y = bx$ passing through $p$. 
It is easy to show that $p$ has the form $(\frac{\sqrt 3}{b+ \sqrt 3}, \frac{\sqrt 3 b}{b+ \sqrt 3}$), where $b$ belongs to the interval $[0, \frac{\sqrt 3}{3}]$. 
Here $b=0$ generates $v_0$, while $b= \frac{\sqrt 3}{3}$ generates $m$.

Our $q$ is the intersection of the segment $v_1v_2$ with a straight line $x=cy$, where $c \in [-\frac{\sqrt 3}{3}, \frac{\sqrt 3}{3}]$. 
So $q$ has the form $(\frac{\sqrt 3}{2} c, \frac{\sqrt 3}{2})$.
An easy calculation shows that the directional coefficient of the straight line containing $pq$ is $\sigma = {{b- \sqrt 3} \over {2 - bc - \sqrt 3 c}}$ and that the directional coefficient of the straight line containing $sp$ is $\varsigma = {{3b+ \sqrt 3} \over {2 + bc + \sqrt 3 c}}$. 

Denote by $S$ the strip between the straight lines containing $pq$ and $rs$, and by $S^+$ the narrowest strip parallel to $S$  which contains $P_6$.
Denote by $T$ the strip between the straight lines containing $qr$ and $sp$, and by $T^+$ the narrowest strip parallel to $T$ which contains $P_6$.

We see that $P = S\cap T$. 
Of course, $S^+ \cap T^+$ is the parallelogram with sides parallel to the sides of $P$ which is circumscribed about $P_6$.
Since we are looking only for parallelograms $P \in \mathcal{P}$, the parallelogram $S^+ \cap T^+$ should be a positive homothetic copy of $P$.

Now, our task is to describe such parallelograms $P$.

We find the first coordinate of the intersection of the straight line through $pq$ with $y = \sqrt 3 x$.
It is $x_0 = \frac{\sqrt 3}{2} \cdot \frac{1-c\sigma}{\sqrt 3 - \sigma}$.
Next we evaluate the ratio of the first 
coordinate of $v_1$ to $x_0$ which is $\frac{\sqrt 3}{3} \cdot \frac{\sqrt 3 - \sigma}{1-c\sigma}$.
By Lemma this ratio equals to $\width(S^+)/\width(S)$.
By the substitution of $\sigma$ we get $\width(S^+)/\width(S) = \frac{\sqrt 3}{3} \cdot \frac{3\sqrt 3 - \sqrt 3 bc -3c - b}{2-2bc}$.
In a similar way we obtain that $\width(T^+)/\width(T) = \frac{\sqrt 3}{3} \cdot \frac{\sqrt 3 + \varsigma}{1- c\varsigma} 
= \frac{\sqrt 3}{3} \cdot \frac{ 3\sqrt 3 + 3 b + 3c+ \sqrt 3 bc}{2-2bc}$.

{\ }

\begin{center}

\includegraphics[width=3.8in]{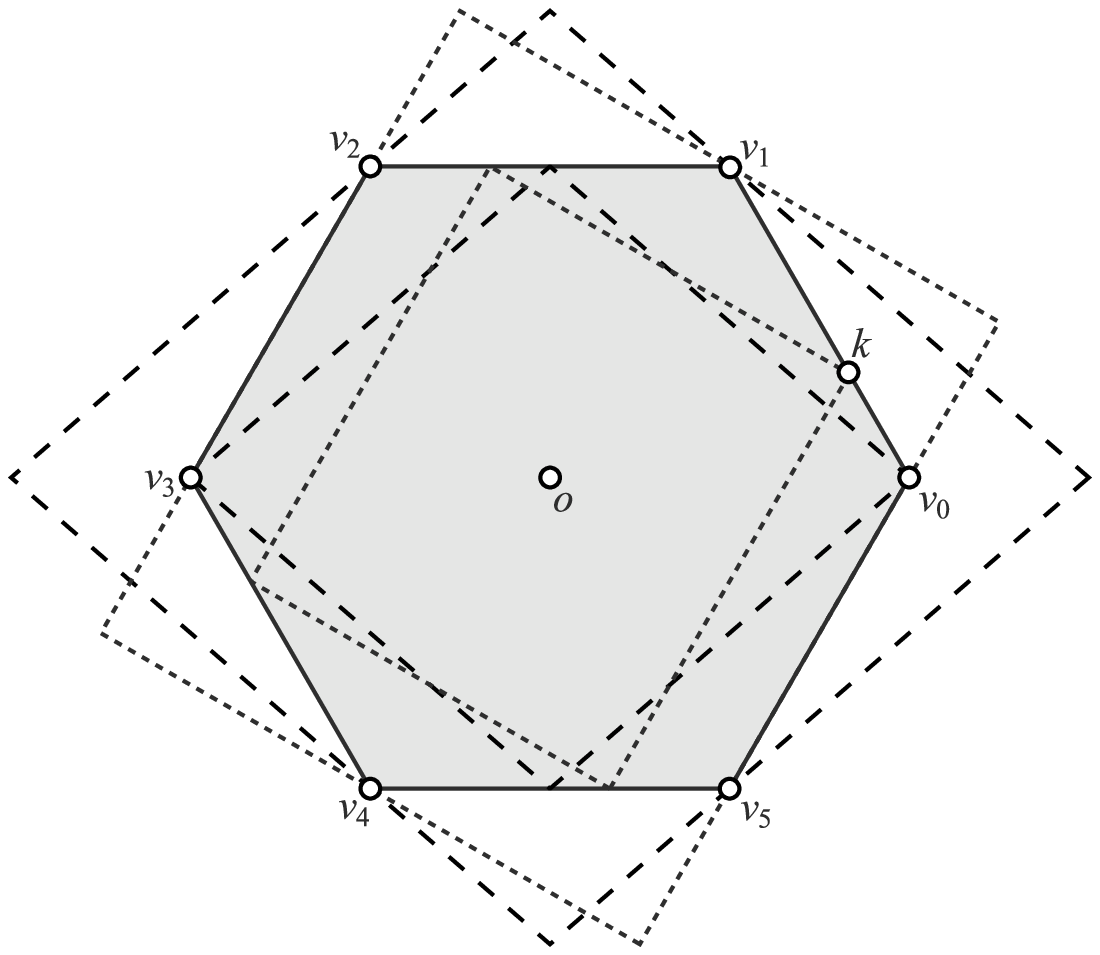} \\ 

\vskip0.39cm
{Fig 2. Two best positions of $P \in \mathcal{P}$ and $h(a)P$ with respect to the hexagon} 

\end{center}

\noindent

Solving the equation $\width(S^+)/\width(S) = \width(T^+)/\width(T)$, we conclude that only for $c = {-2b \over \sqrt 3 b + 3}$ its both sides are equal.
We see that $q \in nv_2$, where $n$ is the midpoint of $v_1v_2$.
Since our $P$ is a function of $b$, we denote it by $P(b)$.
Substituting $c = {-2b \over \sqrt 3 b + 3}$ into $\width(S^+)/\width(S)$ we obtain that the common value of $w(S^+)/w(S)$ and $\width(T^+)/\width(T)$ is $h(b)= {{b^2 + 4\sqrt 3 b + 9} \over {4b^2 +2\sqrt 3 b +6}}$, where $h(b)$ stands in place of $h(P(b))$. 
This ends our task.

Every $P(b)$ is inscribed in $P_6$ and every $h(b)P(b)$ is circumscribed about $P_6$ (see Fig. 1).
The vertices of $h(b)P(b)$ being the images of vertices $p, q, r, s$ of $P(b)$ are denoted by $p', q', t', u'$, respectively. 
Of course, the straight line containing the side $p'q'$ of $h(b)P(b)$ supports $P_6$ at $v_1$.

We are considering here only the situation when the straight line containing the side $s'p'$ of $h(b)P(b)$ supports $P_6$ at $v_5$.
This holds true if and only if the directional 
coefficient $\varsigma$ of the straight line containing $s'p'$ is at most 
the directional coefficient of the straight line containing $v_5v_0$, so when it is at most $\sqrt 3$.
Recall that $\varsigma = {{3b+ \sqrt 3} \over {2 + bc + \sqrt 3 c}}$, which after substituting 
$c = {-2b \over \sqrt 3 b + 3}$ gives $\varsigma= {{3\sqrt 3 b +3} \over {2(\sqrt 3 - b)}}$.
Solving ${{3\sqrt 3 b +3} \over {2(\sqrt 3 - b)}} = \sqrt 3$ we obtain $b= \frac{\sqrt 3}{5}$.

So the straight line containing the side $h(b)\hskip 0.03cm ur$ of $h(b)P(b)$ supports $P_6$ at $v_5$ if and only if $b$ is in the interval $[0, \frac{\sqrt 3}{5}]$.
In other words, if and only if $r\in v_0k$, where $k$ is the point of intersection of $y= \frac{\sqrt 3}{5}x$ with $v_0v_1$.

The derivative of the function $h(b)$ is $\frac{-7\sqrt 3 b^2 - 30 b +3\sqrt 3}{4(b^2 +\sqrt 3 b +6)^2}$. 
An evaluation shows that this derivative is $0$ if and only if $b= {1\over 14}(-10\sqrt 3 \pm \sqrt{384})$.

Only $b= {1\over 14}(-10\sqrt 3 +\sqrt{384}) \approx 0.1625$ in the interval $[0, \frac{\sqrt 3}{5}]$.
The value of $h(b)$ for this $b$ is approximately $1.5224$, so over $\frac{3}{2}$.
Since $h(0) = \frac{3}{2}$ and $h({\sqrt 3\over 5}) = \frac{3}{2}$, we conclude that the global minimum of $h(b)$ in the interval $[0, \frac{\sqrt 3}{5}]$ is $\frac{3}{2}$ and that it is attained only for $b=0$ and $b= \frac{\sqrt 3}{5}$.
In Fig. 2 we see the two pairs $P(b), \frac{3}{2}P(b)$ for these two values of $b$.

It remains to explain what happens for $b \in [\frac{\sqrt 3}{5}, \frac{\sqrt 3}{3}]$, so when $p \in km$.
Observe that then the boundary of the smallest positive homothetic copy of $P$ containing $P_6$ touches it at points $v_1, v_3, v_4, v_0$.
Hence this situation is symmetric to the preceding one with respect to the straight line through $o$ perpendicular to $v_5v_0$. 
Considering $b$ only in the interval $[0, \frac{\sqrt 3}{3}]$ is sufficient since we have in mind 
rotations of $P(b)$ by $60^\circ$ and $120^\circ$ and the axial symmetries with respect to the lines containing $v_0v_3, v_1v_4, v_2v_5$. 
\end{proof}

Let us add that when $b$ changes from $0$ to $\frac{\sqrt 3}{5}$ in the paragraph before the last of the proof, the point $p$ changes from $v_0= (0, 1)$ to $k= (\frac{5}{6}, \frac{\sqrt 3}{6})$ beating $\frac{1}{3}$ of the unit. 
Simultaneously $q$ changes from $(0, \frac{\sqrt 3}{2})$ to $(-\frac{1}{6}$, $\frac{\sqrt 3}{2})$ beating $\frac{1}{6}$ of the unit.

\vskip0.1cm
\noindent
{\bf Remark.}
The only two positions of $P \subset P_6$ such that $P_6 \subset \frac{3}{2}P$ (besides their rotations by $60^\circ$ and $120^\circ$ and axial symmetries with respect to the lines containing $v_0v_3, v_1v_4, v_2v_5$) are the parallelograms with vertices
$(1,0)$, $(0, \frac{\sqrt 3}{2})$, $(-1,0)$, $(0,-\frac{\sqrt 3}{2})$ and with vertices
$(\frac{5}{6}, \frac{\sqrt 3}{6})$, $(-\frac{\sqrt 3}{6}$, $\frac{\sqrt 3}{2})$, $(-\frac{5}{6}$, $-\frac{\sqrt 3}{6})$, $(\frac{\sqrt 3}{6}$, $-\frac{\sqrt 3}{2})$. 
Again see Fig. 2.

\section{Banach-Mazur distance between the parallelogram and the affine-regular even-gons}

In this section we consider the BM-distances from $P_4$ to the regular even-gons with more than six vertices.
In Theorem 2 we find these distances to $P_{8j}$ and $P_{8j+4}$, and we present the estimates from above from $P_4$ to $P_{8j+2}$and $P_{8j+6}$. 
Next we conjecture that the values of these two upper estimates are just the BM-distances from $P_4$ to $P_{8j+2}$ and $P_{8j+6}$.

\begin{thm}\label{I-IV}
We have 

{\rm(I)} $\delta_{BM} (P_4 , P_{8j}) = \sqrt 2$,

{\rm (II)} $\delta_{BM} (P_4 , P_{8j+2}) \leq  \frac{1}{2} \sec \frac{2j}{8j+2}\pi  + \cos \frac{2j}{8j+2}\pi $,

{\rm (III)} $\delta_{BM}(P_4, P_{8j+4}) = \sqrt 2 \cos\frac{1}{8j+4}\pi$,

{\rm (IV)}  $\delta_{BM} (P_4 , P_{8j+6}) \leq  \sin \frac{2j+2}{8j+6}\pi \cdot \csc \frac{4j+2}{8j+6}\pi  + \cos \frac{2j+2}{8j+6}\pi $.
\end{thm}

\begin{proof}
The inequalities showing that the left sides are at most the right sides result from positions of the inscribed parallelograms whose vertices are the intersections of the coordinate axes with the boundaries of $P_{8j}$, $P_{8j+2}$,  $P_{8j+4}$, $P_{8j+6}$, respectively.
Evaluating the ratio of the smallest homothetic copy of such an inscribed parallelogram which contains our polygon from amongst $P_{8j}, P_{8j+2},  P_{8j+4}, P_{8j+6}$, we find each value at the right sides of (I)--(IV). 

Let us show the opposite inequalities for (I) and (III). 

\vskip0.1cm
Ad (I). 
To prove that $\delta_{BM} (P_4 , P_{8j}) \geq \sqrt 2$, we have to show that for any parallelogram $P \in \mathcal M^2$ contained in $P_{8j}$ such that $P_{8j} \subset \lambda P$, where $\lambda > 0$, the inequality $\lambda \geq \sqrt 2$ holds true.
By Proposition, in order to show that $\delta_{BM} (P_4, P_{8j}) \geq  \sqrt 2$, it is sufficient to 
consider only any $P = a(P_4)$ inscribed in $P_{8j}$ such that $\lambda P$ is circumscribed about $P_{8j}$, and to show that $\lambda \geq \sqrt 2$.
Just we may disregard the other $P \subset P_{8j}$.
This is realized in the following Parts ($\alpha$) and ($\beta$).

\vskip 0.1cm
($\alpha$) 
{\it If a parallelogram $P \in \mathcal M^2$ is inscribed in $P_{8j}$ and its homothetic image is circumscribed, then $P$ is a square.} 

Take into account a parallelogram $P \in \mathcal M^2$ inscribed in $P_{8j}$, which is not a square. 
Denote its successive vertices by $a, b, c, d$, when we go counterclockwise.

First let us explain that they must be in some $2j$-th sides of $P_{8j}$.
In the opposite case, some of them, say $a$ and $b$ does not fulfill this.
We may assume that $a \in v_0v_1$ and $b \not \in v_jv_{2j+1}$ 
(besides when $a= v_0$ or $b = v_{2j+1}$).
Provide the straight line parallel to $ab$ which supports $P_{8j}$ at a vertex $v_k$ with $1 < k < 2j$.
Provide the straight line parallel to $bc$ which supports $P_{8j}$ at a vertex $v_l$ with $2j < l < 4j$.
Then the distances from $v_k$ to $ab$ and from $v_l$ to $bc$ are different. 
This means that no homothetic image of $abcd$ has a chance to be circumscribed about $P_{8j}$.

Hence consider the situation when succesive vertices of $P$ are at every $2j$-th side of $P_{8j}$.
Clearly, $a, c$ are symmetric with respect to $o$ and $c$ is in different distances from $a$ and $b$.
Say, let $a \in v_0v_1$ and 
$b \in v_{2j}v_{2j+1}$.
From the symmetry we see that $|v_0a| = |v_{4j}c|$. 
Let $b' \in v_2v_3$ fulfill $|v_oa| = |v_{2j}b'|$ and let $d'$ be opposite to $b'$.
Then $ab'cd'$ is a square inscribed in $P_{8j}$ 
Clearly an enlarged homothetic image of it is circumscribed about $P_{8j}$.
Thus the half-lines with origin at $o$ through $v_1$ and $v_{2j+1}$ intersect 
the sides $ab'$ and $b'c$, respectively, at points $w_1^+, w_{2j+1}^-$, respectively, in equal distances from $o$.
These half-lines intersect the sides $ab$ and $cb$ at points $z_1^+, z_{2j+1}^-$, respectively.
Since $b \not = b'$, then $P$ is not a square [[?rhombus]] and thus $|oz_1| \not = |oz_{2j+1}|$.
Since $|ov_1| = |ov_{2j+1}|$, we get  $|z_1^+ z_1^-|/|w_1^+ w_1^-| \not = |z_{2j+1}^+ z_{2j+1}^-|/|w_{2j+1}^+ w_{2j+1}^-| $ where $w_i^-$, $z_i^-$ denote the points symmetric to $w_i^+$, $z_i^+$.
From Lemma we conclude that the strip between the straight lines through $v_1$ and $v_{4j+1}$, parallel to $ab$ has a different width than the strip between the straight lines through $v_{2j+1}$ and $v_{6j+1}$, parallel to $bc$. 
Hence no homothetic image of $abcd$ is circumscribed about $P_{8j}$.

Consequently, $P$ must be a square.

\vskip0.1cm
($\beta$)
{\it If a square $P$ is inscribed in $P_{8j}$ and $\lambda P$ is circumscribed about $P_{8j}$, then $\lambda \geq \sqrt 2$, and $\lambda = \sqrt 2$ if and only if the vertices of $P$ are at every $2j$-th vertex of $P_{8j}$ or at the middle of every $2j$-th side of $P_{8j}$.} 

Having in mind the axial symmetries of $P_{8j}$, we may limit our considerations to the squares $P$ inscribed in $P_{8j}$ whose one vertex denoted by $p$ is in $v_0m$, where $m$ is the midpoint of $v_0v_1$.

Denote by $k$ the directional factor $\frac{\sin (\pi/4j)}{\cos (\pi/4j) -1}$ 
of the straight line containing the side $v_0v_1$.
So the line has equation $y= k(x-1)+1$.
Clearly,  $p$ is in the intersection of the segment $v_0m$  with a ray $y=bx$, 
where $x \geq 0$ and $0 \leq b \leq \tan {{\pi} \over 8j}$.
We omit an easy calculation showing that $p = (\frac{k}{k-b}, \frac{kb}{k-b})$. 
Since the diagonals of $P$ are orthogonal, the successive vertex $q$ of $P$ is the intersection of the rotated by $90^\circ$   ray  $y=-{1\over b}x$, where $x\geq 0$, with the side $v_{2j}v_{2j+1}$.  
We easily establish that 
$q = (\frac{-kb}{k-b}, \frac{k}{k-b})$.

We provide the straight line containing $pq$. 
Its equation is $y= \frac{b-1}{b+1}x + \frac{k(b^2+1)}{(k-b)(b+1)}$.
Consider its intersection point $u$ with the line $y=x$. 
An easy calculation shows that both coordinates of $u$ are equal to $\frac{k}{2}\cdot \frac{b^2+1}{k-b}$.
Consequently, the ratio of the first coordinate of $v_j$ to the first coordinate of $u$ equals to 
$h(b) = \frac{\sqrt 2}{2} :  \big(\frac{k}{2}\cdot \frac{b^2+1}{k-b}\big) = {\sqrt 2} \cdot \frac {k-b}{k(b^2+1)}$.
By Proposition, the homothetic square $h(b)P$ is circumscribed about $P_{8j}$.

Clearly, we consider the function $h(b)$ in the interval $[0, \tan {{\pi} \over 8j}]$.
Observe that $h(0) = \sqrt 2 = h(\tan\frac{\pi}{8j})$.
Our aim is to show that $h(b) \geq \sqrt 2$ for every $b$ from this interval. 
We find the derivative $h'(b)= \frac{\sqrt 2}{k} \cdot \frac{b^2-2bk-1}{(b^2 +1)^2}$.
An easy calculation shows that $h'(b) = 0$ in the interval $[0, \tan {{\pi} \over 8j}]$ if and only if $b= k + \sqrt{k^2+1}$.
We have $h(k+\sqrt{k^2+1}) = -\frac{\sqrt 2}{k} \cdot \frac{\sqrt{k^2 +1}}{(k+\sqrt{k^2+1})^2 +1}$.
Applying the fact that $k<0$ we show that this value is over $\sqrt 2$.
Hence $h(b) \geq \sqrt 2$ for every $b \in [0, \tan {{\pi} \over 8j}]$ with equality only for $b=0$ and $b= \tan {{\pi} \over 8j}$.

Consequently, the thesis of Part ($\beta$) holds true.

We obtain $\delta_{BM} (P_4 , P_{8j}) \geq \sqrt 2$, which confirms the required opposite inequality of (I).

\vskip0.1cm
Ad (III).
The vertices $v_0, v_{2h+1}, v_{4h+2}, v_{6h+3}$ of $P_{8j+4}$
are the vertices of $P_4$.
Observe that every side of $P_4$ is parallel to every $h$-th side of $P_{8j+4}$. 
So the side $v_0v_{2h+1}$ of $P_4$ is parallel to the side $v_{j}v_{j+1}$ of $P_{8j+4}$, and so on. 
An evaluation shows that the ratio of the first coordinates of the points of the intersection of the line $y=(\tan \frac{\pi}{4}) x$ with these two sides (the points are centers of these two sides) is $\cos\frac{1}{8j+4}\pi \cdot \sec\frac{1}{4} \pi$.
Analogous is true for every side of $P_4$ and the corresponding parallel side of $P_{8j+4}$.
Consequently, by Lemma the homothetic copy of $P_{4}$ with ratio $\cos\frac{1}{8j+4}\pi \cdot \sec\frac{1}{4} \pi$ contains $P_{8j+4}$.
Consequently, $\delta_{BM}(P_4, P_{8j+4}) \leq \cos\frac{1}{8j+4}\pi \cdot \sec\frac{1}{4} \pi$. 
\end{proof}

From Part ($\beta$) we see the only positions of $P$ for which $\delta_{BM} (P_4 , P_{8j}) = \sqrt 2$ is realized.
Observe that the only positions of $P$ for which $\delta_{BM}(P_4, P_{8j+4})$ is realized are these with the vertices at every $(2j+1)$-th vertex of $P_{8j+4}$.

Generalizing (III) we observe that $\delta_{BM}(P_n, P_{hn}) = \cos\frac{1}{hn}\pi \cdot \sec\frac{1}{n} \pi$ for every even $n\geq 4$ and every odd $h \geq 3$. 

In connection with (II) and (IV) we conjecture that
 $\delta_{BM} (P_4 , P_{8j+2}) =  \frac{1}{2} \sec \frac{2j}{8j+2}\pi  + \cos \frac{2j}{8j+2}\pi$ and 
 $\delta_{BM} (P_4 , P_{8j+6}) =  \sin \frac{2j+2}{8j+6}\pi \cdot \csc \frac{4j+2}{8j+6}\pi  + \cos \frac{2j+2}{8j+6}\pi$.

We also do not know the BM-distances between $P_4$ and the regular odd-gons. 
Besides some special cases, the task of finding  the distances $\delta_{BM} (P_m, P_n)$ seems to be very complicated.

\end{document}